\documentclass[oneside,reqno,12pt]{amsart}
\usepackage{amsfonts,amssymb,eucal,amsmath}
\newtheorem{theorem}{Theorem}

\newtheorem{lemma}{Lemma}
\newtheorem{proposition}{Proposition}
\theoremstyle{definition}
\newtheorem{definition}{Definition}
\newtheorem{remark}{Remark}

\renewcommand{\AA}{\mathbb{A}}
\newcommand{\RR}{\mathbb{R}}
\newcommand{\QQ}{\mathbb{Q}}
\newcommand{\ZZ}{\mathbb{Z}}
\newcommand{\NN}{\mathbb{N}}
\renewcommand{\P}{\mathcal{P}}

\newcommand{\M}{{\mathcal{M}}}

\newcommand{\Z}{\mathcal{Z}}

\newcommand{\ve}{\varepsilon}
\newcommand{\set}[1]{\left\{#1\right\}}

\renewcommand{\le}{\leqslant}
\renewcommand{\ge}{\geqslant}

\newcommand{\ov}[1]{\overline{#1}}

\begin{document}

\language=1

\hfill{\tt Journal of Number Theory 111, no.1 (2005) 33-56}

\

\title{On approximation of {\it P}-adic numbers by {\it p}-adic algebraic numbers}

\author{V.V.~Beresnevich}
\author{V.I.~Bernik}
\author{E.I.~Kovalevskaya}



\date{\today}

\maketitle


\bigskip

\section{\bf Introduction}

\bigskip

Throughout $p\ge 2$ is a fixed prime number, $\QQ_p$ is the field
of $p$-adic numbers, $|\omega|_p$ is the $p$-adic valuation of
$\omega\in\QQ_p$, $\mu(S)$ is the Haar measure of a measurable set
$S\subset\QQ_p$, $\AA_{p,n}$ is the set of algebraic numbers of
degree $n$ lying in $\QQ_p$, $\AA_{p}$ is the set of all algebraic
numbers, $\QQ_p^*$ is the extension of $\QQ_p$ containing $\AA_p$.
There is a natural extension of $p$-adic valuation from $\QQ_p$ to
$\QQ_p^*$ \cite{Cassels-1986, Lutz-1955}. This valuation will also
be denoted by $|\cdot|_p$. The disc in $\QQ_p$ of radius $r$
centered at $\alpha$ is the set of solutions of the inequality
$|x-\alpha|_p<r$. Throughout, $\RR_{>a}=\set{x\in\RR:x>a}$,
$\RR_+=\RR_{>0}$ and $\Psi:\NN\to\RR_+$ is monotonic.

Given a polynomial $P(x)=a_nx^n+a_{n-1}x^{n-1}+\ldots+a_1
x+a_0\in\ZZ[x]$ with $a_n\not=0$, $\deg P=n$ is the degree of $P$,
$H(P)=\max_{0\le i\le n}|a_i|$ is the usual height of $P$. Also
$H(\alpha)$ will stand for the usual height of $\alpha\in\AA_p$,
{\em i.e.}\/ the height of the minimal polynomial for $\alpha$.
The notation $X\ll Y$ will mean $X=O(Y)$ and the one of $X\asymp
Y$ will stand for $X\ll Y\ll X$.

In 1989 V.~Bernik \cite{Bernik-1989} proved A.~Baker's conjecture
by showing that for almost all $x\in\RR$ the inequality
$|P(x)|<H(P)^{-n+1}\Psi(H(P))$ has only finitely many solutions in
$P\in\ZZ[x]$ with $\deg P\le n$ whenever and the sum
\begin{equation}\label{e:001}
\sum_{h=1}^{\infty}\Psi(h)
\end{equation}
converges. In 1999 V. Beresnevich \cite{Beresnevich-1999c} showed
that in the case of divergence of (\ref{e:001}) this inequality
has infinitely many solutions.

We refer the reader to
\cite{Beresnevich-Bernik-Kleinbock-Margulis-2002,
BernikDodson-1999, Beresnevich-2002a,
BernikKleinbockMargulis-2001a, Sprindzuk-1979a} for further
development of the metric theory of Diophantine approximation. In
this paper we establish a complete analogue of the aforementioned
results for the $p$-adic case.

\begin{theorem}\label{thm1}
Let $\Psi:\RR_+\to\RR_+$ be monotonically decreasing and
$M_n(\Psi)$ be the set of $\omega\in\QQ_p$ such that the
inequality
\begin{equation}\label{e:002}
|P(\omega)|_p<H(P)^{-n}\Psi(H(P))
\end{equation}
has infinitely many solutions in polynomials $P\in\ZZ[x]$, $\deg
P\le n$. Then $\mu(M_n(\Psi))=0$ whenever the sum\/
{\rm(\ref{e:001})}\/ converges and $M_n(\Psi)$ has full Haar
measure whenever the sum\/ {\rm(\ref{e:001})}\/ diverges.
\end{theorem}

The following is a $p$-adic analogue of Theorem~2 in
\cite{Beresnevich-1999c}.

\begin{theorem}\label{thm2}
Let $\Psi:\RR_+\to\RR_+$ be monotonically decreasing and
$\AA_{p,n}(\Psi)$ be the set of $\omega\in\QQ_p$ such that the
inequality
\begin{equation}\label{e:003}
|\omega-\alpha|_p<H(\alpha)^{-n}\Psi(H(\alpha))
\end{equation}
has infinitely many solutions in $\alpha\in\AA_{p,n}$. Then
$\mu(\AA_{p,n}(\Psi))=0$ whenever the sum\/ {\rm(\ref{e:001})}\/
converges and $\AA_{p,n}(\Psi)$ has full Haar measure whenever the
sum\/ {\rm(\ref{e:001})}\/ diverges.
\end{theorem}

\bigskip

\section{\bf Reduction of Theorem~\ref{thm1}}

\bigskip

We are now going to show that the convergence part of
Theorem~\ref{thm1} follows from the following two theorems. Also
we show that the divergence part of Theorem~\ref{thm1} follows
from Theorem~\ref{thm2}.

\begin{proposition}\label{thm3}
Let $\delta,\xi\in\RR_+$, $\xi<1/2$, $Q\in\RR_{>1}$ and $K_0$ be a
finite disc in $\QQ_p$. Given a disc $K\subset K_0$, let
$E_{1}(\delta,Q,K,\xi)$ be the set of $\omega\in K$ such that
there is a non-zero polynomial $P\in\ZZ[x]$, $\deg P\le n$,
$H(P)\le Q$ satisfying the system of inequalities
\begin{equation}\label{e:004}
\left\{
\begin{array}{l}
 |P(\omega)|_p<\delta Q^{-n-1},\\[1ex]
 |P'(\alpha_{\omega,P})|_p\ge H(P)^{-\xi},
\end{array}
\right.
\end{equation}
where $\alpha_{\omega,P}\in\AA_p$ is the root of $P$ nearest to
$\omega$ $($if there are more than one root nearest to $\omega$
then we choose any of them$)$. Then there is a positive constant
$c_1$ such that for any finite disc $K\subset K_0$ there is a
sufficiently large number $Q_0$ such that
$\mu(E_{1}(\delta,Q,K,\xi))\le c_1\delta\mu(K)$ for all $Q\ge Q_0$
and all $\delta>0$.
\end{proposition}

\begin{proposition}\label{thm4}
Let $\xi,C\in\RR_+$, $K_0$ be a finite disc in $\QQ_p$ and let
$E_{2}(\xi,C,K_0)$ be the set of $\omega\in\QQ_p$ such that there
are infinitely many polynomials $P\in\ZZ[x]$, $\deg P\le n$
satisfying the system of inequalities
\begin{equation}\label{e:005}
\left\{
\begin{array}{l}
 |P(\omega)|_p<C H(P)^{-n-1},\\[1ex]
 |P'(\alpha_{\omega,P})|_p<H(P)^{-\xi}.
\end{array}
\right.
\end{equation}
Then $\mu(E_{2}(\xi,C,K_0))=0$.
\end{proposition}

\begin{proof}[Proof of the convergence part of Theorem~\ref{thm1} modulo
Propositions~\ref{thm3} and \ref{thm4}] Let the sum (\ref{e:001})
converges. Then it is readily verified that
\begin{equation}\label{e:006}
  \sum_{t=1}^\infty 2^t\Psi(2^t)<\infty
\end{equation}
and
\begin{equation}\label{e:007}
  \Psi(h)= o(h^{-1})
\end{equation}
as $h\to\infty$. For the proofs of (\ref{e:006}) see Lemma~5 in
\cite{Beresnevich-1999c}. The arguments for (\ref{e:007}) can be
found in the proof of Lemma~4 in \cite{Beresnevich-1999c}.

Fix any positive $\xi<1/2$. By (\ref{e:007}),
$H(P)^{-n}\Psi(H(P))< H(P)^{-n-1}$ for all but finitely many $P$.
Then, by Proposition~\ref{thm4}, to complete the proof of the
convergence part of Theorem~\ref{thm1} it remains to show that for
any finite disc $K$ in $\QQ_p$ the set $E_{1}(\xi,\Psi)$
consisting of $\omega\in\QQ_p$ such that there are infinitely many
polynomials $P\in\ZZ[x]$, $\deg P\le n$ satisfying the system of
inequalities
\begin{equation}\label{e:008}
\left\{
\begin{array}{l}
 |P(\omega)|_p<H(P)^{-n}\Psi(H(P)),\\[1ex]
 |P'(\alpha_{\omega,P})|_p\ge H(P)^{-\xi}
\end{array}
\right.
\end{equation}
has zero measure.

The system (\ref{e:008}) implies
\begin{equation}\label{e:009}
\left\{
\begin{array}{l}
 |P(\omega)|_p<(2^t)^{-n-1}2^t\Psi(2^t),\\[1ex]
 |P'(\alpha_{\omega,P})|_p\ge H(P)^{-\xi},
\end{array}
\right.
\end{equation}
where $t=t(P)$ with $2^t\le H(P)<2^{t+1}$, which means that
$\omega\in E_{1}(2^{n+1}2^t\Psi(2^t),2^{t+1},K,\xi)$. The system
(\ref{e:009}) holds for infinitely many $t$ whenever (\ref{e:008})
holds for infinitely many $P$. Therefore,
$$
E_{1}(\xi,\Psi)\subset\limsup_{t\to\infty}E_{1}(2^{n+1}2^t\Psi(2^t),2^{t+1},K,\xi).
$$
By Proposition~\ref{thm3},
$\mu(E_{1}(2^{n+1}2^t\Psi(2^t),2^{t+1},K,\xi))\ll 2^t\Psi(2^t)$.
Taking into account (\ref{e:006}), the Borel-Cantelli lemma
completes the proof.
\end{proof}

Next, we are going to show that the divergence part of
Theorem~\ref{thm1} is a consequence of Theorem~\ref{thm2}.

\begin{proof}[Proof of the divergence part of
Theorem~\ref{thm1} modulo Theorem~\ref{thm2}] Fix any finite disc
$K$ in $\QQ_p$. Then there is a positive constant $C>0$ such that
$|\omega|_p\le C$ for all $\omega\in K$. Let $\Psi:\RR_+\to\RR_+$
be a given monotonic function such that the sum (\ref{e:001})
diverges. Then the function $\tilde\Psi(h)=|n!|_p C^{1-n}\Psi(h)$
is also monotonic and the sum $\sum_{h=1}^\infty\tilde\Psi(h)$
diverges. By Theorem~\ref{thm2}, for almost every $\omega\in K$
there are infinitely many $\alpha\in\AA_{p,n}$ satisfying
\begin{equation}\label{e:010}
  |\omega-\alpha|_p<H(\alpha)^{-n}\tilde\Psi(H(\alpha)).
\end{equation}
As $\Psi$ decreases, the right hand side of (\ref{e:010}) is
bounded by a constant. Then we can assume that
$|\omega-\alpha|_p\le C$ for the solutions of (\ref{e:010}). Then
$|\alpha|_p=|\alpha-\omega+\omega|_p\le\max\{|\alpha-\omega|_p,
|\omega|_p\}\le C$.

Let $P_\alpha$ denote the minimal polynomial for $\alpha$. Since
$P_\alpha^{(i)}$ is a polynomial with integer coefficients of
degree $n-i$, we have $|P_\alpha^{(i)}(\alpha)|_p\le \max_{0\le
j\le n-i}|\alpha|_p^{j}\le C^{n-i}$. Then
$$
|P_\alpha(\omega)|_p= |\omega-\alpha|_p\left|\sum_{i=1}^n
i!^{-1}P_\alpha^{(i)}(\alpha)(\omega-\alpha)^{i-1}\right|_p \le
$$
$$
\le |\omega-\alpha|_p\cdot\max_{1\le i\le n}
\left|i!^{-1}P_\alpha^{(i)}(\alpha)(\omega-\alpha)^{i-1}\right|_p\le
$$
$$
\le|\omega-\alpha|_p\cdot|n!|_p^{-1}
C^{n-i}C^{i-1}=|n!|_p^{-1}C^{n-1}|\omega-\alpha|_p.
$$
Therefore (\ref{e:010}) implies
\begin{equation}\label{e:011}
|P_\alpha(\omega)|_p<H(\alpha)^{-n}\tilde\Psi(H(\alpha))|n!|_p^{-1}C^{n-1}=H(\alpha)^{-n}\Psi(H(\alpha))
=H(P_\alpha)^{-n}\Psi(H(P_\alpha)).
\end{equation}
Inequality (\ref{e:010}) has infinitely many solutions for
almost all $\omega\in K$ and so has (\ref{e:011}). As $\omega$
is almost every point of $K$, the proof is completed.
\end{proof}

\bigskip

\section{\bf Reduction of Theorem~\ref{thm2}}

\bigskip

\begin{proof}[Proof of the convergence part of Theorem~\ref{thm2}]
Given an $\alpha\in\AA_{p,n}$, let $\chi(\alpha)$ be the set of
$\omega\in\QQ_p$ satisfying (\ref{e:003}). The measure of
$\chi(\alpha)$ is $\ll H(\alpha)^{-n}\Psi(H(\alpha))$. Then
$$
\sum_{\alpha\in\AA_{p,n}} \mu(\chi(\alpha))=
\sum_{h=1}^\infty\sum_{\alpha\in\AA_{p,n},\ H(\alpha)=h}
\mu(\chi(\alpha))\ll
$$
$$
\ll \sum_{h=1}^\infty\sum_{\alpha\in\AA_{p,n},\ H(\alpha)=h}
h^{-n}\Psi(h)\ll \sum_{h=1}^\infty \Psi(h)<\infty.
$$
Here we used the fact that the quantity of algebraic numbers of
height $h$ is $\ll h^n$. The Borel-Cantelli lemma completes the
proof.
\end{proof}

The proof of the divergence part of Theorem~\ref{thm2} will rely
on the regular systems method of \cite{Beresnevich-1999c}. In this
paper we give a generalization of the method for the $p$-adic
case.

\begin{definition}\label{def1}
Let a disc $K_0$ in $\QQ_p$, a countable set of $p$-adic numbers
$\Gamma$ and a function $N:\Gamma\to\RR_+$ be given. The pare
$(\Gamma,N)$ is called {\em a regular system of points in $K_0$}\/
if there is a constant $C>0$ such that for any disc $K\subset K_0$
for any sufficiently large number $T$ there exists a collection
$$
\gamma_1,\dots,\gamma_t\in\Gamma\cap K
$$
satisfying the following conditions
$$
N(\gamma_i)\le T\ \ (1\le i\le t),
$$
$$
|\gamma_i-\gamma_j|_p\ge T^{-1}\ \ (1\le i<j\le t),
$$
$$
t\ge CT\mu(K).
$$
\end{definition}

\begin{proposition}\label{thm5}
Let $(\Gamma,N)$ be a regular system of points in
$K_0\subset\QQ_p$, ${\tilde\Psi}:\RR_+\to\RR_+$ be monotonically
decreasing function such that $\sum_{h=1}^\infty
{\tilde\Psi}(h)=\infty$. Then $\Gamma_{\tilde\Psi}$ has full Haar
measure in $K_0$, where $\Gamma_{\tilde\Psi}$ consists of
$\omega\in K_0$ such that the inequality
\begin{equation}\label{e:012}
|x-\gamma|_p<\tilde\Psi(N(\gamma))
\end{equation}
has infinitely many solutions $\gamma\in\Gamma$.
\end{proposition}

This theorem is proved in \cite{Beresnevich-Kovalevskaya-2003}.
The proof is also straitforward the ideas of the proof of
Theorem~2 in \cite{Beresnevich-1999c}.

\begin{proposition}\label{thm6}
The pare $(\Gamma,N)$ of\/ $\Gamma=\AA_{p,n}$ and
$N(\alpha)=H(\alpha)^{n+1}$ is a regular system of points in any
finite disc $K_0\subset\QQ_p$.
\end{proposition}

\begin{proof}[Proof of the divergence part of Theorem~\ref{thm2} modulo
Propositions~\ref{thm5} and \ref{thm6}] Let $\Psi:\RR_+\to\RR_+$
be a monotonic function and the sum (\ref{e:001}) diverges. Fix
any finite disc $K_0\subset\QQ_p$.

Let $(\Gamma,N)$ be a regular system defined in
Proposition~\ref{thm5} and let $\Psi$ be a monotonic function such
that the sum (\ref{e:001}) diverges. Define a function
$\tilde\Psi$ by setting
$\tilde\Psi(x)=x^{-n/(n+1)}\Psi(x^{1/(n+1)})$. Using the
monotonicy of $\Psi$, we obtain
$$
\sum_{h=1}^\infty\tilde\Psi(h)=
\sum_{t=1}^\infty\sum_{(t-1)^{n+1}<h\le t^{n+1}}\tilde\Psi(h)
\ge\sum_{t=1}^\infty\sum_{(t-1)^{n+1}<h\le t^{n+1}}t^{-n}\Psi(t)=
$$
$$
=\sum_{t=1}^\infty\left(t^{n+1}-(t-1)^{n+1}\right)t^{-n}\Psi(t)\asymp
\sum_{h=1}^\infty \Psi(h)=\infty.
$$
In is obvious that $\tilde\Psi$ is monotonic. Then, by
Proposition~\ref{thm4}, for almost all $\omega\in K_0$ the
inequality
\begin{equation}\label{e:013}
|x-\alpha|_p<\tilde\Psi(N(\alpha))=H(\alpha)^{-n}\Psi(H(\alpha))
\end{equation}
has infinitely many solutions in $\alpha\in\AA_{p,n}$. The proof
is completed.
\end{proof}

\bigskip

\section{\bf Proof of Proposition~\ref{thm3}}

\bigskip

Fix any finite $K\subset K_0$ in $\QQ_p$. Let $\chi(P)$ be the set
of $\omega\in K$ satisfying (\ref{e:004}) and let $\P_n(Q,K)$ be
the set of non-zero polynomials $P$ with integer coefficients,
$\deg P\le n$, $H(P)\le Q$ and with $\chi(P)\not=\emptyset$. We
will use the following

\begin{lemma}\label{L4}
Let $\alpha_{\omega,P}$ is the nearest root of a polynomial $P$ to
$\omega\in\QQ_p$. Then
$$
|\omega -\alpha_{\omega,P}|_p \le |P(\omega)|_p
|P'(\alpha_{\omega,P})|_p^{-1}.
$$
\end{lemma}

For the proof see \cite[p.\,78]{Sprindzuk-1969}.

Given a polynomial $P\in\P_n(Q,K)$, let $\Z_P$ be the set of roots
of $P$. It is clear that $\#\Z_P\le n$. Given an $\alpha\in\Z_P$,
let $\chi(P,\alpha)$ be the subset of $\chi(P)$ consisting of
$\omega$ with
$|\alpha-\omega|_p=\min\set{|\alpha'-\omega|_p:\alpha'\in\Z_P}$.

By Lemma~\ref{L4}, for any $P\in\P_n(Q,K)$ and any $\alpha\in\Z_P$
one has
\begin{equation}\label{e:014}
\mu(\chi(P,\alpha))\ll \delta Q^{-n-1}|P'(\alpha)|_p^{-1}.
\end{equation}
Given a $P\in\P_n(Q,K)$ and an $\alpha\in\Z_P$, define the disc
\begin{equation}\label{e:015}
\overline\chi(P,\alpha)=\left\{\omega\in K:|\omega-\alpha|_p
\le\Big(4Q|P'(\alpha)|_p\Big)^{-1}\right\}.
\end{equation}
It is readily verified that if
$\overline\chi(P,\alpha)\not=\emptyset$ then
$\mu(\overline\chi(P,\alpha))\gg \Big(4Q|P'(\alpha)|_p\Big)^{-1}$.
Using (\ref{e:014}) we get
\begin{equation}\label{e:016}
\mu(\chi(P,\alpha))\ll \delta Q^{-n-1}\mu(\overline\chi(P,\alpha))
\end{equation}
with the implicit constant depending on $p$ only.

Fix any $P\in\P_n(Q,K)$ and an $\alpha\in\Z_P$ such that
$\chi(P,\alpha)\not=\emptyset$. Let
$\omega\in\overline\chi(P,\alpha)$. Then
\begin{equation}\label{e:017}
P(\omega)=P'(\alpha)(\omega-\alpha)+
(\omega-\alpha)^2\left(\sum_{i=2}^nP^{(i)}(\alpha)
(\omega-\alpha)^{i-2}\right).
\end{equation}

By the inequalities $|P'(\alpha)|_p\ge H(P)^{-\xi}$ and $H(P)\le
Q$, we have $|P'(\alpha)|_p^{-1}\le Q^{\xi}$. Then by
(\ref{e:015}), $|\omega-\alpha|_p\le Q^{-1+\xi}$. Next, as
$\omega\in K$ and $K$ is finite, it is readily verified that
$|P^{(i)}(\alpha)|_p\ll1$, where the constant in this inequality
depends on $K$. Then
\begin{equation}\label{e:018}
\left|(\omega-\alpha)^2\left(\sum_{i=2}^nP^{(i)}(\alpha)
(\omega-\alpha)^{i-2}\right)\right|_p\ll Q^{-2+2\xi}.
\end{equation}
By (\ref{e:015}), we have $|P'(\alpha)(\omega-\alpha)|_p\le
(4Q)^{-1}$. Using this inequality, (\ref{e:018}) and $\xi<1/2$,
we conclude that
\begin{equation}\label{e:019}
|P(\omega)|_p\le (4Q)^{-1},\qquad \omega\in\chi(P,\alpha)
\end{equation}
if $Q$ is sufficiently large.

Assume that $P_1,P_2\in\P_n(Q,K)$ satisfy $P_1-P_2\in\ZZ_{\not=0}$
and assume that there is an
$\omega\in\overline\chi(P_1)\cap\overline\chi(P_2)$. Then
$\omega\in\overline\chi(P_1,\alpha)\cap\overline\chi(P_2,\beta)$
for some $\alpha\in\Z_{P_1}$ and $\beta\in\Z_{P_2}$. Then,
(\ref{e:019}), $|P_1(\omega)-P_2(\omega)|_p<(4Q)^{-1}$. On the
other $P_1(\omega)-P_2(\omega)$ is an integer not greater than
$2Q$ in absolute value. Therefore,
$|P_1(\omega)-P_2(\omega)|_p\ge(2Q)^{-1}$ that leads to a
contradiction. Hence there is no such an $\omega$ and
$\overline\chi(P_1)\cap\overline\chi(P_2)=\emptyset$. Therefore
\begin{equation}\label{e:020}
\sum_{P\in\P_n(Q,K,a_n,\dots,a_1)}\mu(\overline\chi(P))\le\mu(K),
\end{equation}
where $\P_n(Q,K,a_n,\dots,a_1)$ is the subset of $\P_n(Q,K)$
consisting of $P$ with fixed coefficients $a_n,\dots,a_1$.

By (\ref{e:016}) and (\ref{e:020}),
$\sum_{P\in\P_n(Q,K,a_n,\dots,a_1)}\mu(\chi(P))\ll \delta
Q^{-n}\mu(K)$. Summing this over all $(a_n,\dots,a_1)\in\ZZ^n$
with coordinates at most $Q$ in absolute value gives
\begin{equation}\label{e:021}
\sum_{P\in\P_n(Q,K)}\mu(\chi(P))\ll \delta \mu(K).
\end{equation}
It is obvious that
\begin{equation}\label{e:022}
E_{1}(\delta,Q,K,\xi)=\bigcup_{P\in\P_n(Q,K)}\chi(P).
\end{equation}
As the Haar measure is subadditive (\ref{e:021}) and (\ref{e:022})
imply the statement of Proposition~\ref{thm3}.

\language=3

\bigskip

\section{\bf Reduction to irreducible primitive leading polynomials in Proposition~\ref{thm4}}

\bigskip

The following lemma shows us that there is no loss of generality
in neglecting reducible polynomials while proving
Proposition~\ref{thm4}.

\begin{lemma}[Lemma 7 in \cite{Bernik-Dickinson-Yuan-1999}]\label{L1}
Let $\delta\in\mathbb{R}_+$ and $E(\delta)$ be the set of
$\omega\in\QQ_p$ such that the inequality
$$
|P(\omega)|_p <H(P)^{-n-\delta}
$$
has infinitely many solutions in reducible polynomials
$P\in\ZZ[x]$, $\deg P\le n$. Then $\mu (E(\delta))=0$.
\end{lemma}

Also, by Sprind\u{z}uk's theorem \cite{Sprindzuk-1969} there is no
loss of generality in assuming that $\deg P=n$. From now on, $\P$
will denote the set of irreducible polynomials $P\in\ZZ[x]$ with
$\deg P=n$.

\medskip

Next, a polynomial $P\in\ZZ[x]$ is called {\em primitive}\/ if the
gcd (greatest common divisor) of its coefficients is 1. To perform
the reduction to primitive polynomials we fix an $\omega$ such
that the system (\ref{e:005}) has infinitely many solutions in
polynomials $P\in\P$ and show that either $\omega$ belongs to a
set of measure zero or (\ref{e:005}) holds for infinitely many
primitive $P\in\P$.

Define $a_P={\rm gcd}(a_n,\ldots ,a_1,a_0)\in\NN$. Given a
$P\in\P$, there is a uniquely defined primitive polynomial $P_1$
({\em i.e.}\/ $a_{P_1}=1$) with $P=a_P P_1$. Then $H(P)=a_P
H(P_1)$. Let $P\in\P$ be a solution of (\ref{e:005}). By
(\ref{e:005}), $P_1$ satisfies the inequalities
\begin{equation}\label{e:023}
\left\{
\begin{array}{l}
|a_P|_p|P_1(\omega)|_p=|P(\omega)|_p\ll H(P)^{-n-1}=(a_P
H(P_1))^{-n-1},\\[1ex]
|a_P|_p|P_1'(\alpha_{\omega,P})|_p=|P'(\alpha_{\omega,P})|_p<H(P)^{-\xi}=
(a_P H(P_1))^{-\xi}.\\
\end{array}
\right.
\end{equation}
As $|a_P|_p^{-1}\le a_P$, (\ref{e:023}) implies
\begin{equation}\label{e:024}
|P_1(\omega)|_p\ll H(P_1)^{-n-1}a_P^{-n},\ \ \ \
|P_1'(\alpha_{\omega,P})|_p< H(P_1)^{-\xi}a_P^{1-\xi}.
\end{equation}
If (\ref{e:024}) takes place only for a finite number of different
polynomials $P_1\in\P$, then there exists one of them such that
(\ref{e:005}) has infinitely many solutions in polynomials $P$
with the same $P_1$. It follows that $\omega$ is a root of $P_1$
and thus belongs to a set of measure zero. Further we assume that
there are infinitely many $P_1$ satisfying (\ref{e:024}).

If $\xi\ge 1$ then the reduction to primitive polynomials is
obvious as $a_P\in\NN$. Let $\xi<1$. Then, if (\ref{e:005}) holds
for infinitely many polynomials $P\in\P$ such that $a_P\ge
H(P_1)^{\xi'}$, where $\xi'=\xi/(2-2\xi)$, then the first
inequality in (\ref{e:024}) implies that $|P_1(\omega)|_p\ll
H(P_1)^{-n-1}a_P^{-n}\le H(P_1)^{-n-1-n\xi'}$ holds for infinitely
many polynomials $P_1\in\P$. By Sprind\u{z}uk's theorem
\cite{Sprindzuk-1969}, the set of those $\omega$ has zero measure.

If (\ref{e:005}) holds for infinitely many polynomials $P\in\P$
such that $a_P<H(P_1)^{\xi'}$ then (\ref{e:024}) implies that the
system of  inequalities
$$
|P_1(\omega)|_p\ll H(P_1)^{-n-1}, \ \ \ \
|P'(\alpha_{\omega,P})|_p< H(P_1)^{-\xi+(1-\xi)\xi'}<H(P_1)^{-\xi
/2}
$$
holds for infinitely many polynomials $P_1$. Thus, we get the
required statement with a smaller $\xi$.

\medskip

A polynomial $P\in\ZZ[x]$ with the leading coefficient $a_n$ will
be called {\em leading}\/ if
\begin{equation}\label{e:025}
a_n=H(P)\qquad\text{ and }\qquad|a_n|_p>p^{-n}.
\end{equation}

Let $\P_n(H)$ be the set of irreducible primitive leading
polynomials $P\in\ZZ[x]$ of degree $n$ with the height $H(P)=H$.
Also define
\begin{equation}\label{e:026}
\P_n=\bigcup_{H=1}^{\infty}\P_n(H).
\end{equation}

Reduction to leading polynomials is completed with the help of

\begin{lemma}\label{L2}
Let $\Omega$ be the set of points $\omega\in\QQ_p$ for which\/
{\rm(\ref{e:005})}\/ has infinitely many solutions in irreducible
primitive polynomials $P\in\ZZ[x]$, $\deg P=n$. Let $\Omega_0$ be
the set of points $\omega\in\QQ_p$ for which\/
{\rm(\ref{e:005})}\/ has infinitely many solutions in polynomials
$P\in\P_n$, where $\P_n$ is defined in\/ {\rm(\ref{e:026})}\/.
If\/ $\Omega$ has positive measure then so has $\Omega_0$ with
probably a different constant $C$ in\/ {\rm(\ref{e:005})}.
\end{lemma}

Proof of this lemma is very much the same as the one of Lemma~10
in \cite{Sprindzuk-1969} and we leave it as an exercise.

Every polynomial $P\in\P_n$ has exactly $n$ roots, which can be
ordered in any way: $\alpha_{P,1},\dots,\alpha_{P,n}$. The set
$E_{2}(\xi,C,K_0)$ can be expressed as a union of subsets
$E_{2,k}(\xi,C,K_0)$ with $1\le k\le n$, where
$E_{2,k}(\xi,C,K_0)$ is defined to consist of $\omega\in K_0$ such
that (\ref{e:005}) holds infinitely often with
$\alpha_{\omega,P}=\alpha_{P,k}$. To prove Proposition~\ref{thm4}
it suffices to show that $E_{2,k}(\xi,C,K_0)$ has zero measure for
every $k$. The consideration of these sets will not depend on $k$.
Therefore we can assume that $k=1$ and omit this index in the
notation of $E_{2,k}(\xi,C,K_0)$. Also whenever there is no risk
of confusion we will write $\alpha_1,\dots,\alpha_n$ for
$\alpha_{P,1},\dots,\alpha_{P,n}$.

\bigskip

\section{\bf Auxiliary statements and classes of polynomials}

\bigskip

\begin{lemma}\label{L3}
Let $\alpha_1,\ldots ,\alpha_n$ be the roots of $P\in\P_n$. Then
$\max\limits_{1\le i\le n}|\alpha_i|_p <p^n$.
\end{lemma}

For the proof see \cite[p.\,85]{Sprindzuk-1969}.

For the roots $\alpha_1,\ldots ,\alpha_n$ of $P$ we define the
sets
$$
S(\alpha_i)=\{\omega\in \QQ_p : |\omega
-\alpha_i|_p=\min\limits_{1\le j\le n} |\omega
-\alpha_j|_p\}\,\,\,\,\,(1\leq i\leq n).
$$

Let $P\in P_n$. As $\alpha_1$ is fixed, we reorder the other roots
of $P$ so that $|\alpha_1-\alpha_2|_p\le
|\alpha_1-\alpha_3|_p\le\ldots\le|\alpha_1-\alpha_n|_p$. We can
assume that there exists a root $\alpha_m$ of $P$ for which
$|\alpha_1-\alpha_m|_p\le 1$ (see \cite[p.\,99]{Sprindzuk-1969}).
Then we have
\begin{equation}\label{e:027}
|\alpha_1-\alpha_2|_p\le
|\alpha_1-\alpha_3|_p\le\ldots\le|\alpha_1-\alpha_m|_p\le 1\le
\ldots\le |\alpha_1-\alpha_n|_p.
\end{equation}

Let $\ve >0$ be sufficiently small, $d>0$ be a large fixed number
and let  $\ve_1 =\ve /d$, $T=[\ve_1^{-1}]+1$. We define real
numbers $\rho_j$ and  integers $l_j$ by the relations
\begin{equation}\label{e:028}
|\alpha_1-\alpha_j|_p=H^{-\rho_j},\,\,\,\,\,\,(l_j -1)/T\le \rho
_j <l_j /T\,\,\,\,\,\,\, (2\le j \le m).
\end{equation}
It follows from  (\ref{e:027}) and  (\ref{e:028}) that $\rho_2 \ge
\rho_3 \ge \ldots \ge \rho_m \ge 0$ and $l_2 \ge l_3 \ge \ldots
\ge l_m \ge 1$. We assume that $\rho_j =0$ and $l_j =0$ if $m<j\le
n$.

Now for every polynomial $P\in\P_n(H)$ we define a vector $\ov
l=(l_2,\ldots,l_n)$ having non-negative components. In
\cite[p.\,99--100]{Sprindzuk-1969} it is shown that the number of
such vectors is finite and depends on $n, p$ and $T$ only. All
polynomials $P\in {\P_n(H)}$ corresponding to the same vector $\ov
l$ are grouped together into a class ${\P_n(H,\ov l)}$. We define
\begin{equation}\label{e:029}
\P_n(\ov l)=\bigcup_{H=1}^\infty \P_n(H,\ov l).
\end{equation}
Let $K_0=\{\omega\in \QQ_p : |\omega|_p <p^n\}$ be the disc of
radius $p^n$ centered at $0$. Define
$$
r_j=r_j(P)=(l_{j+1}+\ldots +l_n)/T\ \ \ (1\le j \le n-1).
$$

\begin{lemma}\label{L5}
Let $\omega\in S(\alpha_1)$ and $P\in {\P_n(H)}$. Then
$$
H^{-r_1}\ll |P'(\alpha_1)|_p \ll H^{-r_1 +(m-1)\ve_1},
$$
$$
|P^{(j)}(\alpha_1)|_p \ll H^{-r_j +(m-j)\ve_1}\ \ \ \mbox{for} \ \
\ \ 2\le j\le m,
$$
$$
|P^{(j)}(\alpha_1)|_p \ll 1\ \ \ \ \ \ \ \ \ \ \ \ \ \ \ \ \
\mbox{for} \ \ \ \ m<j\le n.
$$
\end{lemma}

\begin{proof}
From (\ref{e:025}) we have $p^{-n}<|H|_p\le 1$. Then, on
differentiating the identity $P(\omega)=H(\omega -\alpha_1)\cdots
(\omega -\alpha_n)$ $j$ times $(1\le j\le n)$ and using
(\ref{e:027}), (\ref{e:028}) we get the statement of the lemma.
\end{proof}

\begin{lemma}\label{L6}
Let $\delta\in\RR_+$, $\sigma\in\RR_+$, $n\ge 2$ be a natural
number and $H=H(\delta ,n)$ be a sufficiently large real number.
Further let $P$, $Q$ in $\ZZ[x]$ be two relatively prime
polynomials of degree at most $n$ with $\max(H(P) ,H(Q))\le H$.
Let $K(\alpha, p^{-t})$ be a disc of radius $p^{-t}$ centered at
$\alpha$ where $t$ is defined by the inequalities $p^{-t}\le
H^{-\sigma}<p^{-t+1}$. If there exists a number $\tau>0$ such that
for all $\omega\in K(\alpha, p^{-t})$ one has
$$
\max(|P(\omega)|_p ,|Q(\omega)|_p)<H^{-\tau}
$$
then $\tau +2\max(\tau -\sigma , 0)<2n+\delta $.
\end{lemma}

For the proof see Lemma~5 in \cite{Bernik-Dickinson-Yuan-1999}.

\bigskip

\section{\bf Proof of Proposition~\ref{thm4}}

\bigskip

As in the previous section
$K_0=\{\omega\in\QQ_p:|\omega|_p<p^n\}$.

Let $A(\ov l,\xi)$ be the set of points $\omega\in K_0$ for which
\begin{equation}\label{e:030}
\left\{
\begin{array}{l}
 |P(\omega)|_p<C H(P)^{-n-1},\\[1ex]
 |P'(\alpha_{1})|_p<H(P)^{-\xi}
\end{array}
\right.
\end{equation}
has infinitely many solutions in polynomials $P\in \P_n({\ov l})$,
where $\P_n({\ov l})$ is defined in (\ref{e:029}). It follows from
the previous discussion that to prove Proposition~\ref{thm4} it
suffices to show that $A(\ov l,\xi)$ has zero measure for all
possible vectors $\ov l$.

The following investigation essentially depends on the value of
$r_1 +l_2/T$. According to Lemma \ref{L5} we have
$|P'(\alpha_1)|_p\gg H^{-r_1}$. It follows from this and the
second inequality of (\ref{e:030}) that $H^{-r_1}\le cH^{-\xi}$,
{\em i.e.}
\begin{equation}\label{e:031}
r_1\ge\xi-\ln c/\ln H>\xi/2 \ \ \ {\rm for}\ \ H\ge H_0.
\end{equation}
Further we assume that $r_1$ satisfies (\ref{e:031}). Further we
set $\ve$ to be $\xi/2$.

\begin{lemma}\label{L7}
If \ $r_1+l_2 /T>n$ then the set of points $\omega\in K_0$ for
which the inequality
$$
|P(\omega)|_p<H(P)^{-n-\ve}
$$
holds for infinitely many polynomials $P\in \P_n({\ov l})$ has
zero measure.
\end{lemma}

For the proof see Proposition 3 in \cite[p.\,111]{Sprindzuk-1969}.

The proof of Proposition~\ref{thm4} is divided into 3 cases, each
corresponding to one of the propositions of this section (see
below).

Let $\chi(P)=\{\omega\in K_0\cap S(\alpha_{P,1}) :
|P(\omega)|_p<H^{-n-1}\}$. Thus, we investigate the set of
$\omega$ that belong to infinitely many $\chi(P)$.

\begin{proposition}\label{P1}
If \   $n-1+2n\ve_1<r_1 +l_2 /T$ then $\mu (A(\ov l,\xi))=0$.
\end{proposition}

\begin{proof}
Let $r_1+l_2 /T>n$. Using Lemma \ref{L7} with $\ve <1$ we get $\mu
(A(\ov l,\xi))=0$.

Let $n-1+2n\ve_1<r_1 +l_2 /T\le n$ and $t$ be a sufficiently large
fixed natural number. We define the set
$$
\M_t(\ov l)=\bigcup_{2^t\le H< 2^{t+1}} \P_n (H,\ov l).
$$
We divide the set $K_0$ into the discs of radius $2^{-t\sigma}$,
where $\sigma =n+1-r_1 -\ve_1$.

First, we consider the polynomials $P\in \M_t(\ov l)$ such that
there is one of the introduced discs, say $K$, such that
$\chi(P)\cap K\not=\emptyset$ and $\chi(Q)\cap K=\emptyset$ for
$Q\in \M_t(\ov l)\smallsetminus\{P\}$. The number of the discs and
respectively the number of the polynomials is at most $p^n
2^{t\sigma}$. From Lemmas \ref{L4} and \ref{L5} we get
$$
\mu(\chi(P))\ll |P(\omega)|_p |P'(\alpha_1)|_p^{-1}\ll
2^{-t(n+1-r_1)}
$$
and thus summing the measures of $\chi(P)$ for the polynomials $P$
of this class leads to
$$
\sum_{P} \mu(\chi(P)) \ll 2^{t(n+1-r_1-\ve_1 -n-1+r_1)}=
2^{-t\ve_1}.
$$
The latest gives the convergent series and, by the Borel-Cantelli
lemma, completes the proof in this case.

Now we consider the other type of polynomials. Let $P$ and $Q$ be
different polynomials of $\M_t(\ov l)$ such that $\chi (P)$ and
$\chi (Q)$ intersect the same disc $D$ introduced above. Then
there exist the points $\omega_1$ and $\omega_2$ belonging to $D$
such that
\begin{equation}\label{e:032}
\max(|P(\omega_1)|_p ,|Q(\omega_2)|_p )\ll 2^{-t(n+1)}.
\end{equation}
Let $\alpha_{P,1}$ and $\alpha_{Q,1}$ be the nearest roots of $P$
and $Q$ to $\omega_1$ and $\omega_2$ respectively. By
(\ref{e:032}), Lemmas \ref{L4} and  \ref{L5} we get
$$
\max(|\omega_1-\alpha_{P,1}|_p ,|\omega_2-\alpha_{Q,1}|_p )\ll
2^{-t(n+1-r_1)}.
$$
Hence, according to the definition of the  $\sigma$ we have
$$
|\alpha_{P,1}- \alpha_{Q,1}|_p\le \max (|\alpha_{P,1}-\omega_1
|_p, |\omega_1 -\omega_2 |_p , |\alpha_{Q,1}-\omega_2 |_p)\ll
$$ $$
\ll\max (2^{-t(n+1-r_1)},2^{-t\sigma})=2^{-t\sigma}.
$$
Now we estimate $|\alpha_{P,1}- \alpha_{Q,i}|_p$\ \ \ \ $(2\le i
\le m)$. Since $r_1 +l_2 /T\le n$ it follows that
$$
|\alpha_{P,1}- \alpha_{Q,i}|_p \le \max (|\alpha_{P,1}-
\alpha_{Q,1}|_p , |\alpha_{Q,1}- \alpha_{Q,i}|_p)\ll \max
(2^{-t\sigma},2^{-t\rho_i})\le
$$
$$
\le \max (2^{-t\sigma},2^{-t(l_i -1)/T})\le 2^{-t(l_i /T-\ve_1)}.
$$
Hence
$$
\prod_{i=1}^m |\alpha_{P,1}-\alpha_{Q,i}|_p\ll 2^{-t(\sigma+(l_2
+\ldots +l_m)/T-(m-1)\ve_1)}= 2^{-t(\sigma+r_1-(m-1)\ve_1)}.
$$
Similarly we obtain
$$
\prod_{i=1}^m |\alpha_{P,2}-\alpha_{Q,i}|_p\le \prod_{i=1}^m
\max\left(
|\alpha_{P,2}-\alpha_{P,1}|_p,|\alpha_{P,1}-\alpha_{Q,1}|_p,|\alpha_{Q,1}-\alpha_{Q,i}|_p\right)\le
$$
$$
\le
\max(2^{-t\rho_2},2^{-t\sigma})\prod_{i=2}^m\max(2^{-t\rho_2},2^{-t\sigma},
2^{-t\rho_i})\ll
$$
$$
\ll 2^{-t(l_2 /T -\ve_1)}\prod_{i=2}^m 2^{-t(l_i /T -\ve_1)} =
2^{-t(l_2 /T-\ve_1+(l_2+\ldots+l_m)/T-(m-1)\ve_1)}=
2^{-t(l_2/T+r_1-m\ve_1)}.
$$
Let $R(P,Q)$ be the resultant of $P$ and $Q$, i.e.
$$
|R(P,Q)|_p =|H|_p^{2n}\prod_{1\le i,j\le n}
|\alpha_{P,i}-\alpha_{Q,j}|_p .
$$
By the previous estimates for $i=1,2$ and the trivial estimates
$|\alpha_{P,i}-\alpha_{Q,j}|_p\ll p^n$ for $3\le i\le n$ we get
$$
|R(P,Q)|_p \ll 2^{-t(\sigma +r_1-(m-1)\ve_1 +l_2 /T+r_1
-m\ve_1)}\le 2^{-t(\sigma +2r_1 +l_2 /T-(2n-1)\ve_1)}<
2^{-t(2n+\delta')}
$$
where $\delta'>0$. On the other hand we have $|R(P,Q)|_p \gg
2^{-2nt}$ as $P$ and $Q$ have not common roots. The last
inequalities lead to a contradiction.
\end{proof}

\begin{proposition}\label{P2}
If
\begin{equation}\label{e:033}
2-\ve /2<r_1+l_2 /T\le n-1+2n\ve_1
\end{equation}
then $\mu (A({\ov l,\xi}))=0$.
\end{proposition}

\begin{proof}
Let
\begin{equation}\label{e:034}
\theta =n+1-r_1 -l_2 /T.
\end{equation}
Let $[\theta]$ and $\{\theta\}$ be the integral and the fractional
parts of $\theta$ respectively.

At first we consider the case  $\{\theta\}\ge\ve$. We define
\begin{equation}\label{e:035}
\beta =[\theta]-1+0,2\{\theta\} -0,1\ve ,
\end{equation}
\begin{equation}\label{e:036}
\sigma_1=l_2 /T +0,8\{\theta\} +(m+1)\ve_1 ,
\end{equation}
\begin{equation}\label{e:037}
d=[\theta]-1.
\end{equation}
Fix any sufficiently large integer $H$ and divide the set $K_0$
into the discs of radius $H^{-\sigma_1}$. The number of these
discs is estimated by $\ll H^{\sigma_1}$. We shall say that the
disc $D$ contains the polynomial $P\in\P_n(H,{\ov l})$ and write
$P\prec D$ if there exists a point $\omega_0\in D$ such that
$|P(\omega_0)|_p<H^{-n-1}$.

Let $B_1 (H)$ be the collection of discs $D$ such that
$\#\{P\in\P_n(H,{\ov l}):P\prec D\}\le H^{\beta}$. By Lemmas
\ref{L4} and \ref{L5}, (\ref{e:035}) and (\ref{e:036}) we have
$$
\sum_{P\in B_1(H)}\mu(\chi(P)) \ll
H^{\beta}H^{\sigma_1}H^{-n-1+r_1}=H^{\theta -1+r_1 +l_2 /T-0,1\ve
+(m+1)\ve_1 -n-1}.
$$
From (\ref{e:034}) we get
$$
\sum_{P}\mu(\chi(P))\ll \sum_{H}H^{-1-\ve/20}<\infty.
$$
By Borel-Cantelli lemma the set of those $\omega$, which belong to
$\chi(P)$ for infinitely many $P\in\bigcup_HB_1(H)$, has zero
measure.

Let $B_2(H)$ be the collection of the discs that do not belong to
$B_1(H)$ and thus contain more than $H^{\beta}$ polynomials
$P\in\P_n(H,{\ov l})$. Let $D\in B_2(H)$. We divide the set
$\{P\in\P_n(H,{\ov l}):P\prec D\}$ into classes as follows. Two
polynomials
$$
P_1(x)=Hx^n+a_{n-1}^{(1)}x^{n-1}+\ldots+a_1^{(1)}x+a_0^{(1)},
$$
$$
P_2(x)=Hx^n+a_{n-1}^{(2)}x^{n-1}+\ldots+a_1^{(2)}x+a_0^{(2)}
$$
are in one class if
$$
a_{n-1}^{(1)}=a_{n-1}^{(2)},\ \ \ldots,\ \ a_{n-d}^{(1)}=a_{n-d}^{(2)},
$$
where $d$ is defined in  (\ref{e:037}). It is clear that the
number of different classes is less than $(2H+1)^d$ and the number
of polynomials under consideration is greater than $H^{\beta}$. By
the pigeon-hole principle, there exists a class $M$ which contains
at least $cH^{\beta -d}$ polynomials where $c>0$ is a constant
independent of $H$. The classes containing less than
$cH^{\beta-d}$ polynomials are considered in a similar way as
above, with the Borel-Cantelli arguments.

Further, we denote polynomials from $M$ by $P_1 (x),\ldots
,P_{s+1}(x)$ and consider $s$ new polynomials
$$
R_1 (x)=P_2 (x)-P_1 (x),\ldots ,R_s(x)=P_{s+1}(x)-P_1(x).
$$
By (\ref{e:037}), we get
\begin{equation}\label{e:038}
\deg R_i\le n-d-1=n-[\theta]\ \ \ \ (1\le i\le s).
\end{equation}
Using (\ref{e:034}), the left-hand side of (\ref{e:033}) and the
condition $\{\theta\}\ge \ve$ we obtain
$$
n-d-1=n-[\theta]=n-\theta +\{\theta\} =-1+r_1 +l_2 /T
+\{\theta\}>1+\ve /2>1.
$$
Since $n-[\theta]$ is integer then
\begin{equation}\label{e:039}
n-[\theta]\ge 2.
\end{equation}
Now we estimate the values $|R_i (\omega)|_p$ $(1\le i\le s)$ when
$\omega\in D$. For every polynomial $P_i$ there exists a point
$\omega_{0i}\in D$ such that  $|P_i (\omega_{0i})|_p<H^{-n-1}$.
Let $\alpha_{1i}$ be the root nearest to  $\omega_{0i}$. By Lemmas
\ref{L4} and \ref{L5}, we get $|\omega _{oi}-\alpha_{1i}|_p\ll
H^{-n-1+r_1}$ and
$$
|\omega -\alpha_{1i}|_p\le \max (|\omega -\omega_{0i}|_p,|\omega
_{0i}-\alpha_{1i}|_p)\ll \max (H^{-\sigma_1},H^{-n-1+r_1})
$$
for any $\omega\in D$. It follows from (\ref{e:036}) and the
right-hand side of  (\ref{e:033}) that
$$
\sigma_1\le n-1-r_1+2n\ve_1 +0,8\{\theta\} +(m+1)\ve_1<n+1-r_1.
$$
Therefore $|\omega -\alpha_{1i}|_p\ll H^{-\sigma_1}$. By Lemma
\ref{L5}, we have
$$
|P_i^{(j)}(\alpha_{1i})(\omega-\alpha_{1i})^{j}|_p\ll
H^{-r_j+(m-j)\ve_1-j\sigma_1}\ \ \ \  \mbox {for}\ \ \  1\le j\le
m,
$$
$$
|P_i^{(j)}(\alpha_{1i})(\omega-\alpha_{1i})^{j}|_p\ll
H^{-j\sigma_1}\ \ \ \ \ \ \ \ \ \ \ \ \  \ \ \ \  \mbox {for}\ \ \
m<j\le n.
$$
From (\ref{e:036}), (\ref{e:034}) and the definition of the $r_j$
$(1\le j\le m)$ we get
$$
|P_i'(\alpha_{1i})(\omega-\alpha_{1i})|_p\ll
H^{-(n+1-\theta)-0,8\{\theta\} -2\ve_1},
$$
$$
|P_i^{(j)}(\alpha_{1i})(\omega-\alpha_{1i})^{j}|_p\ll
H^{-(n+1-\theta)-0,8\{\theta\} -(m+1)\ve_1} \ \ \mbox{for}\ \ 2\le
j\le n.
$$
Using Taylor's formula for $P_i (\omega)$ $(1\le i\le s+1)$ in the
disc $|\omega-\alpha_{1i}|_p\ll H^{-\sigma_1}$ and the previous
estimates, we obtain
\begin{equation}\label{e:040}
|R_i (\omega)|_p\ll H^{-(n+1-\theta)-0,8\{\theta\}
-2\ve_1}=H^{-\tau}\ \ \ (1\le i\le s)
\end{equation}
for any $\omega\in D$. There are the following three cases:
\begin{enumerate}
\item[1)]
Suppose that for each $i$  $(1\le i\le s)$,  $R_i (x)=b_i R(x)$
with $b_i\in\ZZ$. Since the $R_i$ are all different so are the
$b_i$. Let $b=\max\limits_{1\le i\le s}|b_i|=|b_1|$, so that
$b>s/2$. As $bH(R) \le 2H$, $s\gg H^{\beta -d}=H^{0,2\{\theta\}
-0,1\ve}$ and $\{\theta\} \ge\ve$, we get
\begin{equation}\label{e:041}
H(R) \ll H^{1-{0,2\{\theta\} +0,1\ve}} \ \ \ \mbox{and}\ \ \ b\gg
H^{0,2\{\theta\} -0,1\ve}.
\end{equation}
Using (\ref{e:040}) and $H({R_1})=bH(R)$ we have
$$
|R_1 (\omega)|_p=|b|_p |R(\omega)|_p\ll
H({R_1})^{-\tau}=H(R)^{-\tau}b^{-\tau}
$$
and
$$
|R(\omega)|_p\ll H(R)^{-\tau}|b|^{-\tau}|b|_p^{-1}\le
H(R)^{-\tau}b^{-\tau+1}.
$$
From this and (\ref{e:041}) we find
\begin{equation}\label{e:042}
|R(\omega)|_p\ll H(R)^{-\lambda},
\end{equation}
where
$$
\lambda =\tau +(\tau -1)(0,2\{\theta\} -0,1\ve)(1-0,2\{\theta\}
+0,1\ve)^{-1}.
$$
By the definition of the $\tau$ in (\ref{e:040}), the condition
$\{\theta\} \ge\ve$,  (\ref{e:038}) and (\ref{e:039}) we get
$\lambda
>n-[\theta]+1\ge \deg R+1$. It follows from (\ref{e:042}) that
$$
|R(\omega)|_p\ll H(R)^{-\deg R-1-\delta'}
$$
for all $\omega\in D$ , where $\delta'>0$. By Sprind\u{z}uk's
theorem \cite[p.\,112]{Sprindzuk-1969}, the set of $\omega$ for
which there are infinitely many polynomials $R$ satisfying the
previous inequality has zero measure.
\item[2)]
Suppose that some of polynomials $R_i$ are reducible. By
(\ref{e:038}) we have (\ref{e:040}) with  $\tau \ge \deg R_i
+\delta$ where $\delta =1-0,2\{\theta\} +\ve_1
>0$. Then Lemma \ref{L1} shows that the set of $\omega$ for which there are
infinitely many such polynomials has zero measure.
\item[3)]
Suppose that all polynomials $R_i$ are irreducible and that at
least two are relatively prime (otherwise use case 1). Then Lemma
\ref{L6} can be used on two of polynomials, $R_1$ and $R_2$, say.
We have $\deg R_i\le n-[\theta]$ $(i=1,2)$. It follows from
(\ref{e:040}), (\ref{e:034}) and (\ref{e:036})  that
$$
\tau =n+1-\theta +0,8\{\theta\} +2\ve_1=r_1 +l_2 /T +0,8\{\theta\}
+2\ve_1 ,
$$
$$
\tau -\sigma_1 =r_1 -(m-1)\ve_1 =(l_2 +\ldots +l_m)/T-(m-1)/T\ge
T^{-1}>0,
$$
$$
\tau +2(\tau -\sigma_1)=3r_1 +l_2 /T+0,8\{\theta\} -2(m-2)\ve_1,
$$
$$
2(n-[\theta])+\delta =-2+2r_1 +2l_2 /T+2\{\theta\} +\delta .
$$
As $r_1 \ge l_2 /T$ then $\tau +2(\tau
-\sigma_1)>2(n-[\theta])+\delta$ if $0<\delta <\ve $. The last
inequality contradicts Lemma \ref{L6}.
\end{enumerate}

In the case of $\{\theta\} <\ve $ we set
$$
\beta =[\theta]-1+\ve,\ \ \ \ \ \sigma_1=l_2
/T+\{\theta\}+(m+1)\ve_1-(1,5+\ve')\ve, \ \  \ \ve'=\ve /(9n+2)\ \
\ \ \ d=[\theta]-1
$$
and apply the same arguments as above.
\end{proof}

\begin{proposition}\label{P3}
If
\begin{equation}\label{e:043}
\ve \le r_1 +l_2 /T\le 2-\ve /2
\end{equation}
then $\mu (A({\ov l,\xi}))=0$.
\end{proposition}

\begin{proof}
All polynomials
$P(\omega)=H\omega^n+a_{n-1}\omega^{n-1}+...+a_1x+a_0\in\P_n(H,{\ov
l})$ corresponding to the same vector $\ov a=(a_{n-1},\ldots
,a_2)$ are grouped together into a class $\P_n(H,{\ov l},\ov a)$.
Let
$$
B(P)=\{\omega\in K_0\cap S(\alpha_1):|\omega -\alpha_1|_p\le
H^{-n-1} |P'(\alpha_1)|_p^{-1}\},
$$
$$
B_1 (P)=\{\omega\in K_0\cap S(\alpha_1):|\omega -\alpha_1|_p\le
H^{-2+\ve'}|P'(\alpha_1)|_p^{-1}\},
$$
where $\ve'=\ve/6$. It is clear that $B(P)\subset B_1 (P)$,
$$
\mu B(P)=c_1 (p)H^{-n-1}|P'(\alpha_1)|_p^{-1},\ \ \ \ \ \mu B_1
(P)=c_2 (p)H^{-2+\ve'}|P'(\alpha_1)|_p^{-1}
$$
and
\begin{equation}\label{e:044}
\mu B(P)=c_3 (p)H^{-n+1-\ve'}\mu B_1 (P),
\end{equation}
where $c_i (p)>0 \ \ \ (i=1,2,3)$ are the constants dependent on
$p$. Now we estimate $|P(\omega)|_p$ when $P\in \P_n(H,{\ov l},\ov
a)$  and $\omega\in B_1 (P)$. It follows from the definition of
$B_1 (P)$ that $|P'(\alpha_1)(\omega -\alpha_1)|_p<H^{-2+\ve'}$.
By the right-hand side of (\ref{e:043}) and the definition of the
$r_j$ $(2\le j \le m)$ we have
$$
jr_1 -r_j =(j-1)r_1 +r_1 -r_j =(j-1)r_1 +(l_2 +\ldots +l_j)/T\le
(j-1)(2-\ve /2).
$$
From this, Lemma \ref{L5} and the definition of $B_1(P)$ we find
$$
|P^{(j)}(\alpha_1)(\omega -\alpha_1)^j |_p<H^{-r_j +(m-j)\ve
_1}H^{-(2-\ve')j+jr_1}\le H^{-(2-\ve')j+(j-1)(2-\ve
/2)+(m-j)/\ve_1}
=
$$
$$
=H^{-2-(j-1)\ve /2+(m-j)\ve_1+\ve'j}\le H^{-2-\delta}
$$
for $2\le j\le m $, where $\delta >0$ if $\ve_1\le \ve/(2n)$. By
the right-hand side of (\ref{e:043}) and the definition of the
$r_1$ we have $r_1<(2-\ve /2)(1-1/j)$. From this, Lemma \ref{L5}
and the definition of $B_1(P)$ we find
$$
|P^{(j)}(\alpha_1)(\omega -\alpha_1)^j |_p\ll |\omega
-\alpha_1|_p^j<H^{-j(2-\ve'-r_1)}<H^{-2-\ve /3}
$$
for $m<j\le n$. By Taylor's formula and the previous estimates we
get
\begin{equation}\label{e:045}
|P(\omega)|_p \ll H^{-2+\ve'}
\end{equation}
for any $\omega\in B_1 (P)$. Further we use essential and
inessential domains introduced by Sprind\u{z}uk
\cite{Sprindzuk-1969}. The disc $B_1 (P)$ is called {\it
inessential} if there exists a polynomial $Q\in\P_n(H,{\ov l},\ov
a)$ such that $\mu (B_1 (P)\cap B_1 (Q))\ge {\frac12}\mu B_1(P)$
and  {\it essential} otherwise.

Let the disc $B_1(P)$ be inessential and $D=B_1 (P)\cap B_1 (Q)$.
Then
$$ \mu D\ge \textstyle\frac{1}{2}\mu B_1(P)=c_4
(p)H^{-2+\ve'}|P'(\alpha_1)|_p^{-1}
$$
where $c_4(p)>0$ is a constant dependent on $p$. By (\ref{e:045})
the difference $R(\omega)=P(\omega)-Q(\omega)=b_1 \omega +b_0$,
where $\max (|b_0|,|b_1|)\le 2H$, satisfies
\begin{equation}\label{e:046}
|R(\omega)|_p =|b_1\omega -b_0|_p\ll H^{-2+\ve'}
\end{equation}
for any $\omega\in B_1(P)$. Note that $b_1 \neq 0$ since if
$b_1=0$, then  $|b_0|_p \ll H^{-2+\ve'}$. It is contradicted to
$|b_0|_p \ge |b_0|^{-1}\gg H^{-1}$. It follows from (\ref{e:046})
that
\begin{equation}\label{e:047}
|\omega -b_0 /b_1|_p\ll H^{-2+\ve'}|b_1|_p^{-1}.
\end{equation}
Let $D_1 =\{\omega\in K_0\cap S(\alpha_1) :\mbox{the inequality
(\ref{e:047}) holds }\}$. Then $D\subseteq D_1$ and  $\mu D_1 =c_5
(p)H^{-2+\ve'}|b_1|_p^{-1}$, where $c_5 (p)>0$ is a constant
dependent on $p$. We have
$$
c_4 (p)H^{-2+\ve'}|P'(\alpha_1)|_p^{-1}\le \mu D\le \mu D_1 \ll
H^{-2+\ve'}|b_1|_p^{-1}.
$$
Hence
\begin{equation}\label{e:048}
|b_1|_p\ll |P'(\alpha_1)|_p.
\end{equation}
From (\ref{e:048}) and Lemma \ref{L5} we get
$$
|b_1|_p\ll |P'(\alpha_1)|_p\ll H^{-r_1 +(m-1)\ve_1}.
$$
Since $r_1 \ge l_2 /T$ the left-hand side of (\ref{e:043}) implies
$r_1\ge \ve /2$. Now we find $|b_1|_p\ll H^{-\ve /3}$ for
$\ve_1\le \ve/(2n)$. It follows from (\ref{e:046}) that
$|b_0|_p\ll H^{-\ve /3}$. Suppose that $s$ is defined by the
inequalities $p^s\le H<p^{s+1}$. We have $H^{\ve /3}\asymp
p^{[s\ve /3]}$ for sufficiently large $H$. Hence $b_1\asymp
p^{[s\ve /3]}b_{11}$ and $b_0\asymp p^{[s\ve /3]}b_{01}$ where
$b_{11}$, $b_{01}$ are integers. We have
\begin{equation}\label{e:049}
b_1\omega +b_0\asymp p^{[s\ve /3]}(b_{11}\omega +b_{01})\qquad
\mbox{with}\qquad \max(|b_{11}|,|b_{01}|)\ll H^{1-\ve /3}.
\end{equation}
Let $R_1 (\omega)=b_{11}\omega +b_{01}$. Then $H(R_1)\ll H^{1-\ve
/3}$. It follows from (\ref{e:046}) and (\ref{e:049}) that
$$
|b_{11}\omega+b_{01}|_p\ll p^{s\ve /3}H^{-2+\ve'}\ll
H^{-2+\ve'+\ve /3}= H(R_1)^{-2-\ve/(6-2\ve)}.
$$
Using Khintchine's theorem in $\QQ_p$
\cite[p.\,94]{Sprindzuk-1969}, we get that the set of $\omega$
belonging to infinitely many discs $B_1 (P)$ has zero measure.

Let the disc $B_1 (P)$ be essential. By the property of $p$-adic
valuation every point $\omega\in K_0$ belong to no more than one
essential disc. Hence
$$
\sum_{P\in\P(H,\ov l,\ov a)} \mu B_1(P)\le p^n.
$$
It follows from (\ref{e:044}) that
$$
\sum_H\sum_{P\in\P(H,\ov l)} \mu B(P)= \sum_H\sum_{\ov a}
\sum_{P\in\P(H,\ov l,\ov a)} \mu B(P)\ll
$$
$$
\ll \sum_H H^{n-2}\sum_{P\in\P(H,\ov l,\ov a)}H^{-n+1-\ve'}\mu
B_1(P)\ll \sum_H H^{-1-\ve'}< \infty.
$$
The Borel-Cantelli lemma completes the proof.
\end{proof}

\bigskip

\section{\bf Proof of Proposition~\ref{thm6}}

\bigskip

First of all we impose some reasonable limitation on the disc
$K_0$ that appear in the statement of Proposition~\ref{thm6}. To
this end we notice the following two facts.

\begin{remark}
Let $\omega_0,\theta_0\in\QQ_p$. It is a simple matter to verify
that if $(\Gamma,N)$ is a regular system in a disc $K_0$ then
$(\tilde\Gamma,\tilde N)$ is regular in $\theta_0K_0+\omega_0$,
where $\tilde\Gamma=\{\delta_0\gamma+\omega_0:\gamma\in\Gamma\}$,
$\tilde N(\delta_0\gamma+\omega_0)=N(\gamma)$ and
$\theta_0K_0+\omega_0=\{\theta_0\omega+\omega_0:\omega\in K_0\}$.
\end{remark}

\begin{remark}
One more observation is that if $c>0$ is a constant and
$(\Gamma,N)$ is a regular system in a disc $K_0$ then
$(\Gamma,cN)$ is also a regular system in $K_0$.
\end{remark}

The proofs are easy and left as exercises. Now we notice that for
any disc $K_0$ in $\QQ_p$ we can choose two numbers
$\omega_0,\theta_0\in\QQ$ such that $\theta_0\ZZ_p+\omega_0=K_0$.
It is clear that the map $\omega\mapsto\theta_0\omega+\omega_0$
sends $\AA_{p,n}$ to itself. Moreover, there is a constant $c_1>0$
such that for any $\alpha\in\ZZ_p\cap\AA_{p,n}$ one has
$H(\theta_0\alpha+\omega_0)\le c_1H(\alpha)$. Hence, if we will
succeed to prove Proposition~\ref{thm6} for the disc $\ZZ_p$ then
in view of the Remarks above it will be proved for $K_0$. Thus
without loss of generality we assume that $K_0=\ZZ_p$.

In the proof of Proposition~\ref{thm6} we will refer to the
following statement known as Hensel's Lemma (see
\cite[p.~134]{BernikDodson-1999}).

\begin{lemma}
Let $P$ be a polynomial with coefficients in $\ZZ_p$, let
$\xi=\xi_0\in\ZZ_p$ and $|P(\xi)|_p<|P'(\xi)|_p^2$. Then as
$n\to\infty$ the sequence
$$
\xi_{n+1}=\xi_n-\frac{P(\xi_n)}{P'(\xi_n)}
$$
tends to some root $\alpha\in\ZZ_p$ of the polynomial $P$ and
$$
|\alpha-\xi|_p\le\frac{|P(\xi)|_p}{|P'(\xi)|_p^2}<1.
$$
\end{lemma}

\begin{proposition}\label{prp}
Let $\delta>0$, $Q\in\RR_{>1}$. Given a disc $K\subset\ZZ_p$, let
\begin{equation}\label{e:050}
E(\delta,Q,K)=\bigcup_{P\in\ZZ[x],\ \deg P\le n,\ H(P)\le
Q}\{\omega\in K: |P(\omega)|_p<\delta Q^{-n-1}\}.
\end{equation}
Then there is a positive constant $c$ such that for any finite
disc $K\subset\ZZ_p$ there is a sufficiently large number $Q_0$
such that $\mu(E(\delta,Q,K))\le c\delta\mu(K)$ for all $Q\ge
Q_0$.
\end{proposition}

\begin{proof}
The set $E(\delta,Q,K)$ can be expressed as follows
$$
E(\delta,Q,K)\subset E_1(\delta,Q,K,1/3)\,{\textstyle\bigcup}\,
E_3(Q,K)\,{\textstyle\bigcup}\,E_4(),
$$
where $E_1(\delta,Q,K,1/3)$ is introduced in
Proposition~\ref{thm3},
$$
E_3(Q,K)=\bigcup_{P\in\ZZ[x],\ \deg P\le n,\ H(P)\ge \log
Q}\chi(P),
$$
$\chi(P)$ is the set of solutions of (\ref{e:005}) lying in $K$
with $\xi=1/3$ and $C=\delta$,
$$
E_4(Q,K)=\bigcup_{P\in\ZZ[x],\ \deg P\le n,\ H(P)\le \log
Q}\{\omega\in K: |P(\omega)|_p<\delta Q^{-n-1}\}.
$$
By Proposition~\ref{thm4},
\begin{equation}\label{e:051}
\mu(E_3(Q,K))\to0\text{ as }Q\to\infty.
\end{equation}
By Proposition~\ref{thm3},
\begin{equation}\label{e:052}
\mu(E_1(\delta,Q,K,1/3))\le c_1\delta\mu(K)\text{ for sufficiently
large }Q.
\end{equation}
Now to estimate $\mu(E_4(Q,K))$ we first estimate the measure of
$\{\omega\in K: |P(\omega)|_p<\delta Q^{-n-1}\}$ for a fixed $P$.
If $\alpha_{\omega,P}$ is the nearest root to $\omega$ then
$|a_n(\omega-\alpha_{\omega,P})^n|_p<Q^{-n-1}$. Since $|a_n|_p\ge
Q^{-1}$, we get $|\omega-\alpha_{\omega,P}|_p<Q^{-1}$. It follows
that
$$
\mu\{\omega\in K: |P(\omega)|_p<\delta Q^{-n-1}\}\ll Q^{-1}.
$$
Hence $\mu(\mu(E_4(Q,K)))\ll (\log Q)^{n+1}Q^{-1}\to0$ as
$Q\to\infty$. Combining this with (\ref{e:051}) and (\ref{e:052})
completes the proof.
\end{proof}

\begin{proof}[Proof of Proposition~\ref{thm6}]

Fix any disc $K\subset\ZZ_p$ and let $Q>0$ be a sufficiently large
number. Let $\omega\in K$. Consider the system
\begin{equation}\label{e:053}
\left\{
\begin{array}{ll}
|P(\omega)|_p<\delta^2 CQ^{-n-1},&P(\omega)=a_n\omega^n+\dots+a_1\omega+a_0,\\[0.5ex]
|a_j|\le \delta^{-1}Q,& j=\overline{0,n},\\[0.5ex]
|a_j|_p\le \delta,& j=\overline{2,n}.
\end{array}
\right.
\end{equation}
By Dirichlet's principle, it easy to show that there is an
absolute constant $C>0$ such that for any $\omega\in K$ the system
(\ref{e:053}) has a non zero solution $P\in\ZZ[x]$. Fix such a
solution $P$.

If $|P'(\omega)|<\delta$, then, by (\ref{e:053}),
$$
\textstyle |a_1|_p=|P'(\omega)-\sum_{k=2}^nka_k\omega^{k-1}|_p\le
\max\{|P'(\omega)|_p,|2a_2\omega^1|_p,\dots,|na_n\omega^{n-1}|_p\}<\delta.
$$
Also, if $Q$ is sufficiently large, then
$$
\textstyle |a_0|_p=|P(\omega)-\sum_{k=1}^na_k\omega^k|_p\le
\max\{|P(\omega)|_p,|a_1\omega^1|_p,\dots,|a_n\omega^n|_p\}<\delta.
$$
Therefore, the coefficients of $P$ have a common multiple $d$ with
$\delta/p\le|d|_p<\delta$. It follows that $d^{-1}\le \delta$.
Define $P_1=P/d\in\ZZ[x]$. Obviously $H(P_1)\le Q$. Also, by
(\ref{e:053}),
$$
|P_1(\omega)|_p=|P(\omega)|_p|d|_p^{-1}\le
|P(\omega)|_p\times\delta^{-1}p<\delta C p Q^{-n-1}.
$$
This implies $\omega\in E(\delta C p,Q,K)$. By
Proposition~\ref{prp}, $\mu(E(\delta C p,Q,K))\le c \delta C p
\mu(K)$ for sufficiently large $Q$. Put $\delta=(2cpC)^{-1}$. Then
$\mu(K\smallsetminus E(\delta C p,Q,K))\ge\frac12\mu(K)$. If now
we take $\omega\in K\smallsetminus E(\delta C p,Q,K)$ then we get
$$
|P'(\omega)|_p\ge \delta.
$$
By Hensel's lemma there is a root $\alpha\in\ZZ_p$ of $P$ such
that $|\omega-\alpha|_p<CQ^{-n-1}$. If $Q$ is sufficiently large
then $\alpha\in K$. The height of this $\alpha$ is $\le
\delta^{-1}Q$.

Let $\alpha_1,\dots,\alpha_t$ be the maximal collection of
algebraic numbers in $K\cap\AA_{p,n}$ satisfying $H(\alpha_j)\le
\delta^{-1}Q$ and
$$
|\alpha_i-\alpha_j|_p\ge Q^{-n-1}\ \ (1\le i<j\le t).
$$
By the maximality of this collection,
$|\omega-\alpha_j|_p<CQ^{-n-1}$ for some $j$. As $\omega$ is
arbitrary point of $E(\delta C p,Q,K)$, we get
$$
E(\delta C p,Q,K)\subset\bigcup_{j=1}^t
\{\omega\in\ZZ_p:|\omega-\alpha_j|_p<CQ^{-n-1}\}.
$$
Next,
$$
\frac12\mu(K)\le\mu(E(\delta C p,Q,K))\ll Q^{-n-1}t,
$$
whence $t\gg Q^{n+1}\mu(K)$. Taking $T=\delta^{-n-1}Q^{n+1}$ one
readily verifies the definition of regular systems. The proof is
completed.

\end{proof}

{\bf Acknowledgements.} The research was supported by Belorussian
Fond of Fundamental research (Project 00-249) and by INTAS
(project 00-429).



\providecommand{\bysame}{\leavevmode\hbox
to3em{\hrulefill}\thinspace}
\providecommand{\MR}{\relax\ifhmode\unskip\space\fi MR }
\providecommand{\MRhref}[2]{%
  \href{http://www.ams.org/mathscinet-getitem?mr=#1}{#2}
} \providecommand{\href}[2]{#2}

\

{\small

\hangindent=9.5ex \hangafter=1 \noindent {\it Address}\/:
Institute of Mathematics,\\
The National Belarus Academy of Sciences,\\
220072, Surganova 11, Minsk, Belarus

\hangindent=9.7ex \hangafter=1 \noindent {\it E-mails}\/: {\tt
vb8@york.ac.uk, {\it formerly}
beresnevich@im.bas-net.by\\
bernik@im.bas-net.by\\
kovalevsk@im.bas-net.by}

}
\end{document}